\documentclass[a4paper]{amsart}
\usepackage{amssymb}
\def\bfit{\bfseries\itshape}

\newcounter{theo}

\def\equat{\refstepcounter{theo}$$~}
\def\endequat{\leqno{\boldsymbol{(\arabic{theo})}}~$$}

\def\SG{{\mathfrak S}}
\def\RC{{\mathcal{R}}}
\def\NC{{\mathcal{N}}}
\def\a{\alpha}
\def\b{\beta}
\def\D{\Delta}
\def\ph{\varphi}
\def\l{\lambda}
\def\lamt{{\tilde{\lambda}}}
\def\ZM{{\mathbb{Z}}}
\def\FM{{\mathbb{F}}}
\def\NM{{\mathbb{N}}}
\def\lamh{{\hat{\lambda}}}
\def\phba{{\bar{\varphi}}}
\def\lexp#1#2{\kern\scriptspace\vphantom{#2}^{#1}\kern-\scriptspace#2}
\def\le{\hspace{0.1em}\mathop{\leqslant}\nolimits\hspace{0.1em}}
\def\ge{\hspace{0.1em}\mathop{\geqslant}\nolimits\hspace{0.1em}}

\def\DS{\displaystyle}
\DeclareMathOperator{\Ind}{{\mathrm{Ind}}}
\DeclareMathOperator{\Rad}{{\mathrm{Rad}}}
\DeclareMathOperator{\Part}{{\mathrm{Part}}}
\DeclareMathOperator{\Comp}{{\mathrm{Comp}}}
\def\ring{\RC_n}
\def\algebra{\bar{\RC}_n}
\def\radical#1{\bar{\RC}_n^{(#1)}}
\def\gfp{{\FM_{\! p}}}
\def\m{{\mu}}
\def\n{{\nu}}
\def\muh{{\hat{\mu}}}
\def\Delba{{\bar{\Delta}}}
\def\Delt{{\tilde{\Delta}}}
\def\injto{\hookrightarrow}
\def\infspe{\preccurlyeq}
\def\mb{{\mathbf{m}}}
\def\jb{{\mathbf{j}}}
\def\nub{{\boldsymbol{\nu}}}

\begin{document}

\title{A note on the Grothendieck ring \\ of the symmetric group}

\author{C\'edric Bonnaf\'e}
\address{Laboratoire de Math\'ematiques de Besan\c{c}on \\ (CNRS - UMR 6623) \\
Universit\'e de Franche-Comt\'e \\ 16 Route de Gray \\
25030 Besan\c{c}on Cedex\\ France}

\email{bonnafe@math.univ-fcomte.fr}



\date{\today}

\begin{abstract} 
Let $p$ be a prime number and let $n$ be a non-zero natural number. 
We compute the descending Loewy series of the algebra 
$\ring/p\ring$, where $\ring$ denotes the ring 
of virtual ordinary characters of the symmetric group $\SG_n$. 
\end{abstract}

\maketitle

\pagestyle{myheadings}

\markboth{\sc C. Bonnaf\'e}{\sc Grothendieck ring of the symmetric group}


Let $p$ be a prime number and let $n$ be a non-zero natural number. Let 
$\gfp$ be the finite field with $p$ elements and let 
$\SG_n$ be the symmetric group of degree $n$. Let $\ring$ denote 
the ring of virtual ordinary characters of the symmetric group $\SG_n$ 
and let $\algebra=\gfp \otimes_\ZM \ring = \ring/p\ring$. 
The aim of this paper is to determine the 
descending Loewy series of the $\gfp$-algebra $\algebra$ (see Theorem A). 
In particular, we deduce that the Loewy length of $\algebra$ is $[n/p]+1$ 
(see Corollary B). Here, if $x$ is a real number, $[x]$ denotes the unique $r \in \ZM$ 
such that $r \le x < r+1$. 

Let us introduce some notation. If $\ph \in \ring$, 
we denote by $\phba$ its image in $\algebra$. The radical of 
$\algebra$ is denoted by $\Rad \algebra$. If $X$ and 
$Y$ are two subspaces of $\algebra$, we denote by $XY$ the subspace 
of $\algebra$ generated by the elements of the form $xy$, 
with $x \in X$ and $y \in Y$. 

\bigskip

\noindent{\bfit Compositions, partitions.} 
A {\it composition} is a finite sequence $\l=(\l_1,\dots,\l_r)$ 
of non-zero natural numbers. We set $|\l|=\l_1+\dots+\l_r$ and 
we say that $\l$ is a {\it composition of $|\l|$}. 
The $\l_i$'s are called the {\it parts} of $\l$. 
If moreover $\l_1 \ge \l_2 \ge \dots \ge \l_r$, we say that 
$\l$ is a {\it partition of} $|\l|$. The set of compositions 
(resp. partitions) of $n$ is denoted by $\Comp(n)$ (resp. $\Part(n)$).
We denote by $\lamh$ the partition of $n$ obtained from $\l$ 
by reordering its parts. So $\Part(n) \subset \Comp(n)$ 
and $\Comp(n) \to \Part(n)$, $\l \mapsto \lamh$ is surjective. 
If $1 \le i \le n$, we denote by $r_i(\l)$ the number of occurences of 
$i$ as a part of $\l$. We set
$$\pi_p(\l)=\sum_{i=1}^n \Bigl[\frac{r_i(\l)}{p}\Bigr].$$
Recall that $\l$ is called {\it $p$-regular} (resp. {\it $p$-singular}) if and only if 
$\pi_p(\l)=0$ (resp. $\pi_p(\l) \ge 1$). 
Note also that $\pi_p(\l) \in \{0,1,2,\dots,[n/p]\}$ and that 
$\pi_p(\lamh)=\pi_p(\l)$. 
Finally, if $i \ge 0$, we set 
$$\Part_i^{(p)}(n)=\{\l \in \Part(n)~|~\pi_p(\l) \ge i\}.$$

\bigskip

\noindent{\bfit Young subgroups.} 
For $1 \le i \le n-1$, let $s_i=(i,i+1) \in \SG_n$. Let $S_n=\{s_1,s_2,\dots,s_{n-1}\}$. 
Then $(\SG_n,S_n)$ is a Coxeter group. We denote by $\ell : \SG_n \to \NM$ 
the associated length function. If $\l=(\l_1,\dots,\l_r) \in \Comp(n)$, 
we set 
$$S_\l=\{s_i~|~\forall~1 \le j \le r,~i \neq \l_1+\dots+\l_j\}.$$
Let $\SG_\l=<S_\l>$. Then $(\SG_\l,S_\l)$ is a Coxeter group: 
it is a standard parabolic subgroup of $\SG_n$ which is canonically 
isomorphic to $\SG_{\l_1} \times \dots \times \SG_{\l_r}$. 
Note that 
\equat\label{conjugaison sn}
\text{\it $\SG_\l$ and $\SG_\m$ are conjugate in $\SG_n$ if and only if $\lamh=\muh$.}
\endequat
We write $\l \subset \m$ if $\SG_\l \subset \SG_\m$ and we write $\l \le \m$ 
if $\SG_\l$ is contained in a subgroup of $\SG_n$ conjugate to $\SG_\m$. 
Then $\subset$ is an order on $\Comp(n)$ and $\le$ is a preorder on 
$\Comp(n)$ which becomes an order when restricted to $\Part(n)$. 

Let $X_\l=\{w \in \SG_n~|~\forall~x \in \SG_\l,~\ell(wx) \ge \ell(w)\}$. 
Then $X_\l$ is a cross-section of $\SG_n/\SG_\l$. Now, let 
$\NC_\l=N_{\SG_n}(\SG_\l)$ and $W(\l)=\NC_\l \cap X_\l$. 
Then $W(\l)$ is a subgroup of $\NC_\l$ and $\NC_\l=W(\l) \ltimes \SG_\l$. 
Note that 
\equat\label{isomorphisme W}
W(\l) \simeq \SG_{r_1(\l)} \times \dots \times \SG_{r_n(\l)}.
\endequat
Recall that, for a finite group $G$, the {\it $p$-rank} of $G$ is the maximal 
rank of an elementary abelian $p$-subgroup of $G$. 
For instance, $[n/p]$ is the $p$-rank of $\SG_n$. So 
\equat\label{elementary}
\text{\it $\pi_p(\l)$ is the $p$-rank of $W(\l)$.}
\endequat

If $\l$, $\m \in \Comp(n)$, we set 
$$X_{\l\m}=(X_\l)^{-1} \cap X_\m.$$
Then $X_{\l\m}$ is a cross-section of $\SG_\l\backslash\SG_n/\SG_\m$. 
Moreover, if $d \in X_{\l\m}$, there exists a unique composition 
$\n$ of $n$ such that $\SG_\l \cap \lexp{d}{\SG_\m}=\SG_\n$. 
This composition will be denoted by $\l \cap\lexp{d}{\m}$ or 
by $\lexp{d}{\m} \cap \l$. 

\bigskip

\noindent{\bfit The ring ${\boldsymbol{\ring}}$.} 
If $\l \in \Comp(n)$, we denote by $1_\l$ the trivial character of $\SG_\l$ 
and we set $\ph_\l=\Ind_{\SG_\l}^{\SG_n} 1_\l$. Then, by \ref{conjugaison sn}, 
we have $\ph_\l=\ph_{\lamh}$. 
We recall the following well-known old result of Frobenius:
\equat\label{base}
\text{\it $(\ph_\l)_{\l \in \Part(n)}$ is a $\ZM$-basis of $\ring$.}
\endequat
Moreover, by the Mackey formula for tensor product of induced characters, 
we have 
\equat\label{produit tensoriel}
\ph_\l\ph_\m=\sum_{d \in X_{\l\m}} \ph_{\l \cap \lexp{d}{\m}}=\sum_{d \in X_{\l\m}} 
\ph_{\widehat{\l \cap \lexp{d}{\m}}}.
\endequat
Let us give another form of \ref{produit tensoriel}. 
If $d \in X_{\l\m}$, we define 
$\D_d : \NC_\l \cap \lexp{d}{\NC_\m} \to \NC_\l \times \NC_\m$, $w \mapsto (w,d^{-1}wd)$. 
Let $\Delba_d : \NC_\l \cap \lexp{d}{\NC_\m} \to W(\l) \times W(\m)$ be the composition 
of $\D_d$ with the canonical projection $\NC_\l \times \NC_\m \to W(\l) \times W(\m)$. 
Then the kernel of $\Delba_d$ is $\SG_{\l \cap \lexp{d}{\m}}$, so $\Delba_d$ 
induces an injective morphism $\Delt_d : W(\l,\m,d) \injto W(\l) \times W(\m)$, 
where $W(\l,\m,d)=(\NC_\l \cap \lexp{d}{\NC_\m})/\SG_{\l \cap \lexp{d}{\m}}$. 
Now, $W(\l) \times W(\m)$ acts on $\SG_\l \backslash \SG_n /\SG_\m$ and 
the stabilizer of $\SG_\l d \SG_\m$ in $W(\l) \times W(\m)$ 
is $\Delt_d(W(\l,\m,d))$. Moreover, if $d'$ is an element of $X_{\l\m}$ 
such that $\SG_\l d \SG_\m$ and $\SG_\l d' \SG_\m$ are in the same 
$(W(\l)\times W(\m))$-orbit, then $\SG_{\l \cap \lexp{d}{\m}}$ and 
$\SG_{\l \cap \lexp{d'}{\m}}$ are conjugate in $\NC_\l$. Therefore,
\equat\label{nouvelle formule}
\ph_\l \ph_\m = \sum_{d \in X_{\l\m}'} \frac{|W(\l)|.|W(\m)|}{|W(\l,\m,d)|} 
\ph_{\l \cap \lexp{d}{\m}},
\endequat
where $X_{\l\m}'$ denotes a cross-section of $\NC_\l \backslash \SG_n /\NC_\m$ 
contained in $X_{\l\m}$. 

\bigskip

\noindent{\bfit The descending Loewy series of ${\boldsymbol{\algebra}}$.} 
We can now state the main results of this paper.

\bigskip

{\sc Theorem A.} 
{\it If $i \ge 0$, we have
$\DS{\bigl(\Rad \algebra\bigr)^i = 
\mathop{\oplus}_{\l \in \Part_i^{(p)}(n)} \gfp \phba_\l}$.}

\bigskip

{\sc Corollary B.}
{\it The Loewy length of $\algebra$ is $[n/p]+1$.}

\bigskip

Corollary B follows immediately from Theorem A. The end of this paper 
is devoted to the proof of Theorem A. 

\bigskip

\noindent{\bfit Proof of Theorem A.} 
Let $\radical{i}=\DS{\mathop{\oplus}_{\l \in \Part_i^{(p)}(n)} \gfp \phba_\l}$. 
Note that 
$$0=\radical{[n/p]+1}\subset\radical{[n/p]} \subset \dots \subset \radical{1} \subset 
\radical{0}=\algebra.$$ 
Let us first prove the following fact:

$$\text{\it If $i$, $j \ge 0$, then $\radical{i}\radical{j} \subset \radical{i+j}$.}
\leqno{(\clubsuit)}$$
\begin{proof}[Proof of $(\clubsuit)$]
Let $\l$ and $\m$ be two compositions of $n$ such that $\pi_p(\l) \ge i$ and 
$\pi_p(\m) \ge j$. Let $d \in X_{\l\m}'$ be such that $p$ does not divide 
$\DS{\frac{|W(\l)|.|W(\m)|}{|W(\l,\m,d)|}}$. By \ref{nouvelle formule}, 
we only need to prove that this implies that $\pi_p(\l \cap \lexp{d}{\m}) \ge i+j$. 
But our assumption on $d$ means that $\Delt_d(W(\l,\m,d))$ contains a Sylow 
$p$-subgroup of $W(\l) \times W(\m)$. 
In particular, the $p$-rank of $W(\l,\m,d)$ is greater than or equal to 
the $p$-rank of $W(\l) \times W(\m)$. By \ref{elementary}, this means 
that the $p$-rank of $W(\l,\m,d)$ is $\ge i+j$. Since 
$W(\l,\m,d)$ is a subgroup of $W(\l \cap \lexp{d}{\m})$, we get that the $p$-rank 
of $W(\l \cap \lexp{d}{\m})$ is $\ge i+j$. In other words, again by \ref{elementary}, 
we have $\pi_p(\l \cap \lexp{d}{\m}) \ge i+j$, as desired. 
\end{proof}

\bigskip

By $(\clubsuit)$, $\radical{i}$ is an ideal of $\algebra$ and, 
if $i \ge 1$, then $\radical{i}$ is a nilpotent ideal of $\algebra$. 
Therefore, $\radical{1} \subset \Rad \algebra$. In fact:

$$\Rad \algebra=\radical{1}.\leqno{(\diamondsuit)}$$
\begin{proof}[Proof of $(\diamondsuit)$]
First, note that $\Rad \algebra$ consists of the nilpotent elements of 
$\algebra$ because $\algebra$ is commutative. Now, let 
$\ph$ be a nilpotent element of $\algebra$. Write 
$\ph=\sum_{\l \in \Part(n)} a_\l \phba_\l$ and let $\l_0 \in \Part(n)$ be maximal 
(for the order $\le$ on $\Part(n)$) such that $a_{\l_0} \neq 0$. 
Then, by \ref{nouvelle formule}, the coefficient of $\ph_{\l_0}$ in 
$\ph^r$ is equal to $a_{\l_0}^r |W(\l_0)|^{r-1}$. 
Therefore, since $\ph$ is nilpotent and $a_{\l_0} \neq 0$, we get that 
$p$ divides $|W(\l_0)|$, so that $\l_0 \in \Part_1^{(p)}(n)$ 
(by \ref{elementary}). 
Consequently, $\ph - a_{\l_0} \phba_{\l_0}$ is nilpotent and we can repeat 
the argument to find finally that $\ph \in \radical{1}$. 
\end{proof}

\bigskip

We shall now establish 
a special case of \ref{produit tensoriel} (or \ref{nouvelle formule}). 
We need some notation. If $\a=(\a_1,\dots,\a_r)$ is a composition of 
$n'$ and $\b=(\b_1,\dots,\b_s)$ is a composition of $n''$, let $\a \sqcup \b$ 
denote the composition of $n'+n''$ equal to $(\a_1,\dots,\a_r,\b_1,\dots,\b_r)$. 
If $1 \le j \le n$ and if $0 \le k \le [n/j]$, we denote 
by $\nub(n,j,k)$ the composition $(n-jk,j,j,\dots,j)$ of $n$, where $j$ is repeated 
$k$ times (if $n=jk$, then the part $n-jk$ is omitted). If $\l \in \Comp(n)$, we set 
$$M(\l)=\{0\} \cup \{1 \le j \le n~|~\text{$p$ does not divide $r_j(\l)$}\},$$
$$\mb(\l)=\max M(\l),$$
$$J(\l)=\{0\} \cup \{1 \le j \le n~|~r_j(\l) \ge p\},$$
$$\jb(\l)=\min J(\l)$$
$$\jb\mb(\l)=(\jb(\l),\mb(\l)).\leqno{\text{and}}$$
Let $I=\{0,1,\dots,[n/p]\}$. 
Then $\jb\mb(\l) \in I \times I$. Let us now introduce an order $\infspe$ on 
$I \times I$. If $(j,m)$, $(j',m')$ are two elements of $I \times I$, we 
write $(j,m) \infspe  (j',m')$ if one of the following two conditions is 
satisfied:

\begin{centerline}
{\begin{tabular}{l}
(a) $j < j'$. \\
(b) $j=j'$ and $m \ge m'$.
\end{tabular}}
\end{centerline}

Now, let $i \ge 1$ and let $\l \in \Part_{i+1}^{(p)}(n)$. Let 
$(j,m)=\jb\mb(\l)$. Then $\l=\widehat{\a \sqcup \n_0}$, 
where $\a$ is a partition of $n-m-jp$ 
and $\n_0=\nub(m+jp,j,p)$. Let $\lamt=\a \sqcup (m+jp)$. 
Then $\pi_p(\lamt)=i$ (indeed, $r_{m+jp}(\lamt)=1+r_{m+jp}(\a)=1+r_{m+jp}(\l)$ and, 
by the maximality of $m$, we have that $p$ divides $r_{m+jp}(\l)$) and 

$$\phba_{\nub(n,j,p)} \phba_\lamt \in \phba_\l + 
\Bigl(\mathop{\oplus}_{\substack{\m \in \Part_{i+1}^{(p)}(n) \\ 
\jb\mb(\m) \prec (j,m)}} 
\gfp \phba_\m\Bigr).\leqno{(\heartsuit)}$$
\begin{proof}[Proof of $(\heartsuit)$]
Let $\n=\nub(n,j,p)$. 
Since $\pi_p(\l)=i+1 \ge 2$, there exists $j' \in \{1,2,\dots,[n/p]\}$ such that 
$r_{j'}(\l) \ge p$. Then $j' \ge j$ (by definition of $j$), so $n \ge 2pj$. 
In particular, $n-pj > j$, so $W(\n) \simeq \SG_p$. Now, if $m' > m$, then 
$r_{m'}(\a)=r_{m'}(\l)$, so 
$$\forall~m'>m,~r_{m'}(\a) \equiv 0 \mod p.\leqno{(\heartsuit')}$$
Also
$$\forall~l \neq m+jp,~r_l(\lamt)=r_l(\a)\leqno{(\heartsuit'')}$$
and
$$r_{m+jp}(\lamt)=r_{m+jp}(\a)+1.\leqno{(\heartsuit''')}$$
Now, keep the notation of \ref{nouvelle formule}. We may, and we will, 
assume that $1 \in X_{\n\lamt}'$. First, note that $\n \cap \lamt=\a \sqcup \n_0$ 
and that the image of $\Delt_1$ in  
$W(\n) \times W(\lamt)$ is equal to $W(\n) \times W(\a)$. 
But, by $(\heartsuit')$, $(\heartsuit'')$ and $(\heartsuit''')$, 
the index of $W(\a)$ in $W(\lamt)$ is $\equiv 1 \mod p$. 
Thus, by \ref{nouvelle formule}, we have 
$$\phba_\n \phba_\lamt = \phba_\l + 
\sum_{d \in X_{\n\lamt}'-\{1\}} \frac{|W(\n)|.|W(\lamt)|}{|W(\n,\lamt,d)|} 
\phba_{\n \cap \lexp{d}{\lamt}}.$$
Now, let $d$ be an element of $X_{\n\lamt}$ such that $p$ does not divide 
$\DS{\frac{|W(\n)|.|W(\lamt)|}{|W(\n,\lamt,d)|}}=x_d$ and such that 
$\jb\mb(\n \cap \lexp{d}{\lamt}) \succcurlyeq \jb\mb(\l)$. It is sufficient to 
show that $d \in \NC_\n\NC_\lamt$. 
Write $\a=(\a_1,\dots,\a_r)$. Then 
$$\n \cap \lexp{d}{\lamt}=(n_1,\dots,n_r,n_0) \sqcup j^{(1)} \sqcup \dots \sqcup j^{(p)},$$ 
where $n_k \ge 0$ and $j^{(l)}$ is a composition of $j$ with at most $r+1$ parts. 
Since $p$ does not divide $x_d$, the image of 
$\NC_\n \cap \lexp{d}{\NC_\lamt}$ in $W(\n)$ contains a Sylow $p$-subgroup 
of $W(\n)\simeq \SG_p$. Let 
$w \in \NC_\n \cap \lexp{d}{\NC_\lamt}$ be such that its image 
in $W(\n)$ is an element of order $p$. Then there exists $\sigma \in \SG_\n$ such 
that $w\sigma$ normalizes $\SG_{\n \cap \lexp{d}{\lamt}}$. In particular, 
$\widehat{j^{(1)}}=\dots=\widehat{j^{(p)}}$. So, 
if $j^{(1)} \neq (j)$, then $\jb(\n \cap \lexp{d}{\lamt}) < \jb(\l)$, 
which contradicts our hypothesis. 
So $j^{(1)}=\dots=j^{(p)}=(j)$. Therefore, 
$$\lamt \cap \lexp{d^{-1}}{\n} = \nub(\a_1,j,k_1) \sqcup \dots \sqcup \nub(\a_r,j,k_r) 
\sqcup \nub(m+jp,j,k_0),$$
where $0 \le k_i \le p$ and $\sum_{i=0}^r k_i=p$. Note that 
$(n_1,\dots,n_r,n_0)=(\a_1-k_1j,\dots,\a_r-k_rj,\a_0-k_0j)$ where, 
for simplification, we denote $\a_0=m+jp$. Also, 
$\jb(\lamt \cap\lexp{d^{-1}}{\n}) \le j$ and, since 
$\jb\mb(\lamt \cap \lexp{d^{-1}}{\n}) \succcurlyeq (j,m)$, we have that 
$\mb(\lamt \cap \lexp{d^{-1}}{\n}) \le m$. 
Recall that $d^{-1}wd \in \NC_\lamt$. So two cases may occur:

$\bullet$ Assume that there exists a sequence $0 \le i_1 < \dots < i_p \le r$ 
such that $0 \neq k_{i_1}=\dots=k_{i_p}$ ($=1$) and such that $\a_{i_1}=\dots=\a_{i_p}$. 
So $r_l(\lamt \cap \lexp{d^{-1}}{\n})\equiv r_l(\lamt) \mod p$ for every $l \ge 1$. 
In particular, $r_{m+jp}(\lamt \cap \lexp{d^{-1}}{\n}) \equiv 1 + r_{m+jp}(\a) 
\equiv 1 \mod p$ by $(\heartsuit')$ and $(\heartsuit''')$. Thus, 
$\mb(\lamt \cap \lexp{d^{-1}}{\n}) \ge m+jp > m$, which contradicts 
our hypothesis.

$\bullet$ So there exists a unique $i \in \{0,1,\dots,r\}$ such that $k_i=p$. 
Consequently, $k_{i'}=0$ if $i' \neq i$. If $\a_i > m+jp$, then 
$r_{\a_i}(\lamt \cap \lexp{d^{-1}}{\n})=r_{\a_i}(\lamt)-1=r_{\a_i}(\a)-1$ (by 
$(\heartsuit'')$), so $p$ 
does not divide $r_{\a_i}(\lamt \cap \lexp{d^{-1}}{\n})$ (by $(\heartsuit')$), 
which implies that $\mb(\lamt \cap \lexp{d^{-1}}{\n}) \ge \a_i > m$, 
contrarily to our hypothesis. If $\a_i < m+jp$, then 
$r_{m+jp}(\lamt \cap \lexp{d^{-1}}{\n})=r_{m+jp}(\lamt)=r_{m+jp}(\a)+1$ 
(by $(\heartsuit''')$), so $p$ does not divide 
$r_{m+jp}(\lamt \cap \lexp{d^{-1}}{\n})$ (by $(\heartsuit')$), contrarily to 
our hypothesis. This shows that $\a_i=m+jp$. In other words, 
$d \in \NC_\lamt\NC_\n$, as desired. 
\end{proof}

\bigskip

By $(\diamondsuit)$, Theorem A follows immediately from the next result: 
if $i \ge 0$, then 

\bigskip

$$\radical{1} \radical{i} = \radical{i+1}.\leqno{(\spadesuit)}$$
\begin{proof}[Proof of $(\spadesuit)$]
We may assume that $i \ge 1$. By $(\clubsuit)$, we have 
$\radical{1} \radical{i} \subset \radical{i+1}$. 
So we only need to prove that, if $\l \in \Part_{i+1}^{(p)}(n)$ then 
$\phba_\l \in \radical{1}\radical{i}$. But this follows from $(\heartsuit)$ and 
an easy induction on $\jb\mb(\l) \in I \times I$ (for the order $\infspe$).
\end{proof}

\end{document}